\documentclass{elsart}               

\usepackage{graphicx}          
\usepackage{amssymb}    
\usepackage{amsmath}                           
\begin{document}

\begin{frontmatter}

\title{On an intrinsic formulation of time-variant Port Hamiltonian systems}


\author[JKU]{Markus Sch\"{o}berl}
\ead{markus.schoeberl@jku.at}
\author[JKU]{Kurt Schlacher}
\ead{kurt.schlacher@jku.at}
\address[JKU]{Johannes Kepler University Linz, Institute of Automatic Control and Control Systems Technology,
     Austria}


\begin{abstract}                          
In this contribution we present an intrinsic description of time-variant Port Hamiltonian systems as they
appear in modeling and control theory. This formulation is based on the splitting of the state bundle and the use of appropriate covariant derivatives, which guarantees that
the structure of the equations is invariant with respect to time-variant coordinate transformations. 
In particular, we will interpret our covariant system representation in the context of control theoretic problems. 
Typical examples are time-variant error systems related to  trajectory tracking problems which allow for a Hamiltonian formulation. 
Furthermore we will analyze the concept of collocation and the balancing/interaction of power flows in an intrinsic fashion.

\end{abstract}

\end{frontmatter}

\section{Introduction}

Hamiltonian systems are the object of analysis for a long period and
they have been investigated from many different points of view and
in many different scientific areas. In the last two decades, in mathematical
physics especially field theoretic aspects of Hamiltonian systems
without control input are of importance, see \cite{Giachetta,Gotay,Kanat}.
In field theory the use of bundles to distinguish dependent and independent
coordinates is commonly used and since time-variant lumped parameter
systems can be seen as a special case of field theory with only one
independent coordinate the use of bundles also applies to time-variant
systems where the fibration is accomplished with respect to the time-coordinate.
Beside field theory, also in classical mechanics, especially in the
time-invariant setting, the geometric interpretation of the Hamiltonian
picture is well established, see for example \cite{Abraham} for many
details concerning this subject. \\
From the control theoretic point of view, the class of Port Hamiltonian
systems are a well-analyzed class, see for example \cite{Ortega,vanderSchaft}
and references therein, where both the theoretical point of view and,
of course, the physical applications play a prominent role. Roughly
speaking, the main idea of many passivity based control approaches
is to maintain the Hamiltonian structure of the system by feedback
since this structure has some pleasing properties concerning the stability
proof also in the nonlinear scenario. 

In the literature most of these approaches for the lumped parameter
scenario concerning control theoretic aspects present system analysis,
modeling and control for time-invariant systems, whereas the time-variant
case is analyzed very rarely. We believe that the main difficulty
in the time-variant scenario is the fact that the geometric picture
of the equations changes considerably. In the time-invariant setting
the role of the time is solely to be the curve parameter, which is
not true in the time-variant scenario. Contributions which treat the
time-variant case especially with regard to control theoretic problems
are for example \cite{FujimotoCSL,FujimotoAut} where the authors
consider, what they call generalized Hamiltonian systems and canonical
transformations which might be time-dependent. 

Two important applications where time-variant systems arise quite
naturally based on a time-variant change of coordinates should be
mentioned at this stage: Firstly, the introduction of displacement
coordinates with respect to a system trajectory as it arises for instance
when the analysis of the tracking error is the objective, see for
example \cite{FujimotoAut}. And secondly, in mechanics/robotics floating/accelerated
frames of reference are commonly used with respect to an inertial
one. 

The main contribution of this paper are that: (i) an intrinsic definition
of time-variant (Port) Hamiltonian systems is given based on a covariant
derivative induced by a connection; (ii) this intrinsic description
is analyzed in a differential geometric way; (iii) for the system
class of time-variant (controlled) Hamiltonian mechanics, a covariant
version of the power balance relation including collocation is developed;
(iv) for the special case where (beside a possible feed-forward) the
connection can be expressed as an additive Hamiltonian the results
of \cite{FujimotoCSL} are recovered. 

It is worth mentioning that in our opinion a time-variant (Port) Hamiltonian
system has to be introduced using covariant derivatives, which differs
significantly from the definition in \cite{FujimotoCSL,FujimotoAut}.
We identify 'covariant' with the fact that system properties do not
depend on the chosen coordinate chart, i.e., we formulate systems
in an intrinsic way. The key idea is the use of a connection which
induces a covariant derivative, see \cite{Giachetta}. Partially,
results in this paper have been presented preliminarily in \cite{SchoeberlPAMM2006,SchoeberlNolcos}.

\section{The Time-Invariant Case}

This introductory section is a reminder of time-invariant Port Hamiltonian
systems \cite{MaschkeOrtegaSchaft2000TAC,Ortega,vanderSchaft} including
also the arising matching conditions when state transformations and
affine input transformations are considered. It serves as a basis
for the generalization to the time-variant case and is also used to
introduce the differential geometric language which is then extensively
exploited in the time-variant scenario. The notation is similar to
the one in \cite{Giachetta}, where the interested reader can find
much more details about this geometric machinery. 

To keep the formulas short and readable we will use tensor notation
and especially Einsteins convention on sums where we will not indicate
the range of the used indices when they are clear from the context.
We use the standard symbol $\otimes$ for the tensor product, $\mathrm{d}$
is the exterior derivative, $\rfloor$ the natural contraction between
tensor fields and $\circ$ denotes the composition of maps. By $\partial_{A}^{B}$
are meant the partial derivatives with respect to coordinates with
the indices $_{B}^{A}$. \\
To study the time-invariant case of Port Hamiltonian systems in a
geometric fashion we introduce the state manifold $\mathcal{X}$ equipped
with coordinates $(x^{\alpha})$, where $\alpha=1,\ldots,\dim(\mathcal{X})$
and we consider diffeomorphisms (in the sequel also called transition
functions) of the type $\bar{x}=\varphi(x)$ where $\bar{x}$ denotes
the states in the transformed coordinate system. Standard differential
geometric constructions, see \cite{Abraham,Giachetta,NijmeijervanderSchaft,Saunders}
lead to the tangent bundle $\mathcal{T(X)}$ and the cotangent bundle
$\mathcal{T}^{\ast}\mathcal{(X)}$, which possess the induced coordinates
$(x^{\alpha},\dot{x}^{\alpha})$ and $(x^{\alpha},\dot{x}_{\alpha})$
with respect to the holonomic bases $\partial_{\alpha}$ and $\mathrm{d}x^{\alpha}.$
Typical elements of $\mathcal{T}(\mathcal{X})$ (vector fields) and
$\mathcal{T}^{*}(\mathcal{X})$ (\textit{1-forms}) read in local coordinates
as $w=\dot{x}^{\alpha}(x)\partial_{\alpha}$ and $\omega=\dot{x}_{\alpha}(x)\mathrm{d}x^{\alpha}$,
respectively. To introduce in- and outputs we consider the vector
bundle $\mathcal{U}\rightarrow\mathcal{X}$ with the coordinates $(x^{\alpha},u^{i})$
for $\mathcal{U}$ and the base $e_{i}$ for the fibres where $i=1,\ldots,\dim(\mathcal{U_{F}})$,
where $\mathcal{U_{F}}$ denotes the fibres of the input bundle (vector
spaces) as well as the dual output vector bundle $\mathcal{Y}\rightarrow\mathcal{X}$
possessing the coordinates $(x^{\alpha},y_{i})$ and the fibre base
$e^{i}$. Greek indices will correspond to the components of the coordinates
of the state manifold and induced structures. Latin indices correspond
to the components of the input and the output variables (fibres of
the dual bundles $\mathcal{U}\rightarrow\mathcal{X}$ and $\mathcal{Y}\rightarrow\mathcal{X}$).
Let us consider the maps $J,R:\mathcal{T}^{\ast}\mathcal{(X)\rightarrow T(X)}$
which are contravariant tensors that are given by the local coordinate
expressions 
\begin{equation}
J=J^{\alpha\beta}\partial_{\alpha}\otimes\partial_{\beta}\,,\,\,\, R=R^{\alpha\beta}\partial_{\alpha}\otimes\partial_{\beta}\label{eq:JandR}
\end{equation}
with $J^{\alpha\beta},R^{\alpha\beta}\in\mathcal{C}^{\infty}(\mathcal{X})$
where $J$ is skew-symmetric, i.e. $J^{\alpha\beta}=-J^{\beta\alpha}$
and $R$ is symmetric $R^{\alpha\beta}=R^{\beta\alpha}$ and positive-semidefinite.
Furthermore we introduce the bundle map $G:\mathcal{U}\rightarrow\mathcal{T(X)}$
which is a tensor that has the local coordinate expression $G=G_{i}^{\alpha}e^{i}\otimes\partial_{\alpha}$
with $G_{i}^{\alpha}\in C^{\infty}(\mathcal{X})$. Having the maps
$J,R$ and $G$ at our disposal a time-invariant Port Hamiltonian
system (with dissipation), see \cite{MaschkeOrtegaSchaft2000TAC,Ortega,vanderSchaft}
can be constructed as
\begin{equation}
\begin{array}{ccl}
\dot{x} & = & (J-R)\rfloor\mathrm{d}H+G\rfloor u\\
y & = & G^{*}\rfloor\mathrm{d}H
\end{array}\label{eq:HamTimeInvStandardCF}
\end{equation}
where the function $H\in C^{\infty}(\mathcal{X})$ denotes the Hamiltonian
and $G^{*}:\mathcal{T}^{*}(\mathcal{X})\rightarrow\mathcal{Y}$ the
adjoint (dual) map of $G$. The local coordinate expression of (\ref{eq:HamTimeInvStandardCF})
reads as
\begin{equation}
\begin{array}{ccl}
\dot{x}^{\alpha} & = & \left(J^{\alpha\beta}-R^{\alpha\beta}\right)\partial_{\beta}H+G_{i}^{\alpha}u^{i}\\
y_{i} & = & G_{i}^{\alpha}\partial_{\alpha}H.
\end{array}\label{eq:HamTimeInvStandard}
\end{equation}
We want to analyze structure preserving transformations for the system
(\ref{eq:HamTimeInvStandardCF}). To allow for affine input transformations
we can replace the input bundle by an affine one $\mathcal{Z\rightarrow X}$
(with underlying vector bundle $\mathcal{U}\rightarrow\mathcal{X}$),
for the geometric properties of affine bundles see for example \cite{Giachetta}
and references therein. The transition functions for the vector bundle
and the affine bundle read as
\begin{eqnarray}
 &  & \bar{u}=Mu\,,\,\,\,\bar{u}^{\bar{j}}=M_{i}^{\bar{j}}u^{i}\label{eq:UtransLin}\\
 &  & \bar{u}=Mu+g\,,\,\,\,\bar{u}^{\bar{j}}=M_{i}^{\bar{j}}u^{i}+g^{\bar{j}}\label{eq:UTransAffine}
\end{eqnarray}
with $M_{i}^{\bar{j}},g^{\bar{j}}\in C^{\infty}(\mathcal{X})$ where
$\bar{u}$ denotes the transformed input coordinates and we restrict
ourselves to regular transformations (i.e. $M$ is invertible). The
geometric representation of the system leads to the observation that
the structure of (\ref{eq:HamTimeInvStandardCF}) is preserved by
a diffeomorphism of the type $\bar{x}=\varphi(x)$ together with (\ref{eq:UtransLin}).
The case of an affine input bundle is more challenging since the preservation
of the structure demands to solve a partial differential equation.
See also \cite{Cheng} in this context, where the problem of general
feedback equivalence of nonlinear systems to Port Hamiltonian systems
is discussed and so called \textit{matching conditions} appear. 
\begin{lem}
\label{lem:affineInput}Consider the system (\ref{eq:HamTimeInvStandardCF})
together with the diffeomorphism $\bar{x}=\varphi(x)$ and (\ref{eq:UTransAffine}).
The structure of (\ref{eq:HamTimeInvStandardCF}) is preserved if
an only if we can find a solution $\breve{H}\in C^{\infty}(\bar{\mathcal{X}})$
of the partial differential equations 
\begin{equation}
(\bar{J}^{\bar{\alpha}\bar{\beta}}-R^{\bar{\alpha}\bar{\beta}})\partial_{\bar{\beta}}\breve{H}-\left(\partial_{\alpha}\varphi^{\bar{\alpha}}G_{i}^{\alpha}\hat{M}_{\bar{j}}^{i}g^{\bar{j}}\right)\circ\hat{\varphi}=0.\label{eq:PDEaff}
\end{equation}
Here $\bar{J}^{\bar{\alpha}\bar{\beta}}$ and $\bar{R}^{\bar{\alpha}\bar{\beta}}$
are the components of the transformed tensors (\ref{eq:JandR}) with
respect to \textup{$\bar{x}=\varphi(x)$}. The inverse maps are denoted
by $x=\hat{\varphi}(\bar{x})$ and $M_{i}^{\bar{j}}\hat{M}_{\bar{j}}^{k}=\delta_{i}^{k}$
where $\delta$ is the Kronecker delta.\end{lem}
\begin{rem}
The partial differential equations (\ref{eq:PDEaff}) are written
in the coordinates $\bar{x}$ but it is readily observed that it can
be formulated in the original coordinates $x$, as well.
\end{rem}
The proof of this Lemma is a straightforward calculation in local
coordinates. If in Lemma \ref{lem:affineInput} a solution for $\breve{H}$
can be obtained, then the following Corollary is an immediate consequence.
\begin{cor}
Suppose (\ref{eq:PDEaff}) is met, then the system (\ref{eq:HamTimeInvStandard})
in the new coordinates reads as
\[
\begin{array}{ccl}
\dot{\bar{x}}^{\bar{\alpha}} & = & (\bar{J}^{\bar{\alpha}\bar{\beta}}-R^{\bar{\alpha}\bar{\beta}})\partial_{\bar{\beta}}(\bar{H}-\breve{H})+\bar{G}_{\bar{i}}^{\bar{\alpha}}\bar{u}^{\bar{i}}\\
\bar{y}_{\bar{i}} & = & \bar{G}_{\bar{i}}^{\bar{\alpha}}\partial_{\bar{\alpha}}(\bar{H}-\breve{H})
\end{array}
\]
with $\bar{G}_{\bar{j}}^{\bar{\alpha}}=\left(\partial_{\alpha}\varphi^{\bar{\alpha}}G_{i}^{\alpha}\hat{M}_{\bar{j}}^{i}\right)\circ\hat{\varphi}$
and $\bar{H}=H\circ\hat{\varphi}$. The output 
\begin{equation}
\bar{y}_{\bar{j}}=\hat{M}_{\bar{j}}^{i}(y_{i}-(\partial_{\alpha}\breve{H})G_{i}^{\alpha}).
\end{equation}
and the Hamiltonian are transformed affine and in general $u^{i}y_{i}\neq\bar{u}^{\bar{j}}\bar{y}_{\bar{j}}$
is met.\end{cor}
\begin{exmp}
Let us consider the shifting of a nonzero (but constant) equilibrium
point $(x_{e},u_{e})$ of the system (\ref{eq:HamTimeInvStandardCF}).
This leads to affine relations of the form \textup{$\bar{x}=x-x_{e}$
and $\bar{u}=u-u_{e}$. }For a different interpretation concerning
Port Hamiltonian systems with nonzero equilibrium consider \cite{OrtegaEnergy}
whereas in \cite{MaschkeOrtegaSchaft2000TAC} the construction of
Lyapunov functions is treated in this case.
\end{exmp}

\section{The Time-Variant Case}

Time-variant systems arise when the time coordinate is treated as
an additional coordinate in contrast to the case where it is solely
a curve parameter. This has the consequence for (Port) Hamiltonian
systems that the maps $J,R,G$ and the Hamiltonian $H$ may depend
on the time coordinate and that also time-dependent coordinate changes
have to be taken into account. Concerning time-dependent changes of
coordinates in mathematical physics, especially in mechanics a covariant
treatment of the equations requires to use frames of references formulated
on bundles, see \cite{Giachetta}. 
\begin{rem}
A special case of a time-variant coordinate transformation is the
introduction of displacement coordinates with respect to a system
trajectory of the (time-invariant) system (\ref{eq:HamTimeInvStandard}).
\end{rem}
We will apply intrinsic concepts to model  Hamiltonian systems. To
motivate how these covariant concepts arise let us consider a time-variant
system modeled on an extended state manifold given by the direct product
$\mathcal{B\times X}$ that includes the time coordinate $t^{0}$
for $\mathcal{B}$ (the index 0 will always correspond to the time
coordinate in the sequel) as it is analyzed for instance in \cite{FujimotoCSL}
and by a slight abuse of notation we have
\[
\begin{array}{ccl}
\partial_{0}x^{\alpha}(t) & = & J^{\alpha\beta}(t,x(t))\partial_{\beta}H(t,x(t))+G_{i}^{\alpha}(t,x(t))u^{i}(t)\\
y_{i}(t) & = & G_{i}^{\alpha}(t,x(t))\partial_{\alpha}H(t,x(t)).
\end{array}
\]
The crucial point is now how to interpret the time derivative $\partial_{0}x^{\alpha}(t)$
either using the structure of jet spaces (jet bundles) which are affine
or using tangent structures. Let us consider a transformation of the
form $\bar{x}^{\bar{\alpha}}=\varphi^{\bar{\alpha}}(x^{\beta},t^{0})$
with $\bar{t}^{0}=t^{0}$. We have 
\begin{equation}
\dot{\bar{x}}^{\bar{\alpha}}=\partial_{0}\varphi^{\bar{\alpha}}\dot{t}^{0}+\partial_{\beta}\varphi^{\bar{\alpha}}\dot{x}^{\beta}\label{eq:TransTXDot}
\end{equation}
and consequently the vector field $\partial_{0}$$ $ is mapped to
the vector field $\partial_{0}+\partial_{0}\varphi^{\bar{\alpha}}\partial_{\bar{\alpha}}$.
This follows from $\dot{t}^{0}=1$ together with (\ref{eq:TransTXDot})
and the main reason for that fact is the observation that a trivial
product bundle structure of the form $\mathcal{B\times X\rightarrow B}$
is not preserved by time-variant transformations and therefore not
the adequate choice for time-variant systems. To overcome this problem
the system has to be formulated on a bundle $\mathcal{E}\rightarrow\mathcal{B}$.
The choice of the trivial one $\mathcal{B\times X\rightarrow B}$
corresponds to a specific trivialization $\mathcal{E}\cong\mathcal{B}\times\mathcal{X}$,
i.e. to the choice of a certain reference frame. An intrinsic approach
demands a formulation on the bundle $\mathcal{E}\rightarrow\mathcal{B}$
and the choice of a connection, which can be seen as a reference frame.

\subsection{Geometric Setting}

To obtain an appropriate geometric picture for the time-variant case
we consider the bundle $\mathcal{E}\rightarrow\mathcal{B}$, where
we use coordinates $(t^{0},x^{\alpha})$ for $\mathcal{E}$ and obtain
the following geometric structures. The tangent bundle $\mathcal{T(E)}$
with coordinates $(t^{0},x^{\alpha},\dot{t}^{0},\dot{x}^{\alpha})$,
the cotangent bundle $\mathcal{T}^{\ast}\mathcal{(E)}$ equipped with
coordinates $(t^{0},x^{\alpha},\dot{t}_{0},\dot{x}_{\alpha}),$ as
well as the vertical bundle $\mathcal{V(E)}$ with coordinates $(t^{0},x^{\alpha},\dot{x}^{\alpha})$
and $\mathcal{V}^{\ast}\mathcal{(E)}$ with $(t^{0},x^{\alpha},\dot{x}_{\alpha}),$
where $\mathcal{V(E)}$ possess the induced bases $\partial_{\alpha}$
and since $\mathcal{V}^{\ast}\mathcal{(E)}$ does not possess a canonical
base without the choice of a connection it is denoted as $\mathrm{\tilde{d}}x^{\alpha}$
at this stage, see \cite{Giachetta}. \\
To be able to deal with time-derivatives of sections $s:\mathcal{B}\rightarrow\mathcal{E}$,
i.e $x=s(t)$ we introduce the first jet manifold $\mathcal{J}^{1}\mathcal{(E)}$
with the adapted coordinates $(t^{0},x^{\alpha},x_{0}^{\alpha})$,
see \cite{Giachetta,Saunders}. The coordinates $x_{0}^{\alpha}$
are often called derivative coordinates and the transition functions
with respect to the bundle morphism $\bar{x}^{\bar{\alpha}}=\varphi^{\bar{\alpha}}\left(t,x\right),$
$\bar{t}^{0}=t^{0}$ (we do not consider time-reparametrization) read
as $\bar{x}_{0}^{\bar{\alpha}}=\partial_{0}\varphi^{\bar{\alpha}}+\partial_{\beta}\varphi^{\bar{\alpha}}x_{0}^{\beta}$.
\\
The key object for a covariant system representation will be a connection
together with the associated covariant differential. The main philosophy
behind a connection is the fact that $\mathcal{T}(\mathcal{E})$ possesses
a canonical subbundle, namely the vertical tangent bundle $\mathcal{V}(\mathcal{E})$
but there is no canonical horizontal complement such that $\mathcal{T(E)}=\mathcal{V(E)}\oplus\mathcal{H(E)}$
holds, unless one specifies a horizontal%
\footnote{It should be noted that the concept of vertical and horizontal parts
of several objects will become important in the sequel and we will
indicate this using $\mathcal{V}$ for vertical and $\mathcal{H}$
for horizontal. However, $ $$\mathcal{H}$ should not be confused
with the Hamiltonian $H$.%
}subbundle and this is exactly what a connection does. By duality one
derives analogously a decomposition of the cotangent bundle as $\mathcal{T^{*}(E)}=\mathcal{V^{*}(E)}\oplus\mathcal{H^{*}(E)}$,
where in this case $\mathcal{H}^{*}(\mathcal{E})$ is canonically
given and the connection defines $\mathcal{V}^{*}(\mathcal{E})$.

\subsubsection{Connection }

Given a bundle $\mathcal{E}\rightarrow\mathcal{B}$ whose fibration
induces the vertical tangent bundle $\mathcal{V}(\mathcal{E})$ it
is the desire to obtain a splitting of the form $\mathcal{T}(\mathcal{E})=\mathcal{V}(\mathcal{E})\oplus\mathcal{H}(\mathcal{E})$.
This can be achieved by a connection $\Gamma$ which in local coordinates
can be represented as a tensor of the form 
\begin{equation}
\Gamma=\mathrm{d}t^{0}\otimes\left(\partial_{0}+\Gamma_{0}^{\alpha}\partial_{\alpha}\right)\label{eq:Conn}
\end{equation}
with the connection coefficients $\Gamma_{0}^{\alpha}\in\mathcal{C}^{\infty}\mathcal{(E)}$,
see \cite{Giachetta,Saunders}. The transition functions for $\Gamma$
when no time-reparametrization (i.e. $\bar{t}=t$) takes place read
as 
\begin{equation}
\bar{\Gamma}_{0}^{\bar{\alpha}}=\partial_{0}\varphi^{\bar{\alpha}}+\Gamma_{0}^{\alpha}\partial_{\alpha}\varphi^{\bar{\alpha}}\label{eq:connTransLaw}
\end{equation}
 where $\bar{\Gamma}_{0}^{\bar{\alpha}}$ are the connection coefficients
in the transformed coordinate system.

This connection (\ref{eq:Conn}) can be used to define a covariant
derivative and to split the tangent bundle $\mathcal{T}(\mathcal{E})\rightarrow\mathcal{E}$,
i.e. $\mathcal{T(E)}=\mathcal{V(E)}\oplus\mathcal{H(E)}$ as stated
above, where $\mathcal{H}(\mathcal{E})$ denotes the horizontal subbundle.
We have the following coordinate representations for the splitting
for a typical element of $\mathcal{T}(\mathcal{E})$ and $\mathcal{T}^{*}(\mathcal{E})$,
respectively 
\begin{equation}
\begin{array}{ccc}
\dot{t}^{0}\partial_{0}+\dot{x}^{\alpha}\partial_{\alpha} & = & \dot{t}^{0}w_{0}^{\mathcal{H}}+\left(\dot{x}^{\alpha}-\dot{t}^{0}\Gamma_{0}^{\alpha}\right)\partial_{\alpha}\\
\dot{t}_{0}\mathrm{d}t^{0}+\dot{x}_{\alpha}\mathrm{d}x^{\alpha} & = & \left(\dot{t}_{0}+\dot{x}_{\alpha}\Gamma_{0}^{\alpha}\right)\mathrm{d}t^{0}+\dot{x}_{\alpha}\omega_{\mathcal{V}}^{\alpha}
\end{array}\label{eq:splitvec1a}
\end{equation}
with the vector field $w_{0}^{\mathcal{H}}$ and the \textit{1-form}
$\omega_{\mathcal{V}}^{\alpha}$ 
\begin{equation}
w_{0}^{\mathcal{H}}=\partial_{0}+\Gamma_{0}^{\alpha}\partial_{\alpha},\ \omega_{\mathcal{V}}^{\alpha}=\mathrm{d}x^{\alpha}-\Gamma_{0}^{\alpha}\mathrm{d}t^{0}\label{eq:whwv}
\end{equation}
such that $ $$w_{0}^{\mathcal{H}}$ and $\omega_{\mathcal{V}}^{\alpha}$
qualify as bases for $\mathcal{H}(\mathcal{E})$ and $\mathcal{V}^{*}(\mathcal{E})$,
respectively.
\begin{rem}
It is readily observed that a time-variant transformation $\bar{x}^{\bar{\alpha}}=\varphi^{\bar{\alpha}}\left(t,x\right)$
converts a trivial connection $\Gamma_{0}^{\alpha}=0$ in one coordinate
chart to a non-trivial one in the transformed coordinates $\bar{\Gamma}_{0}^{\bar{\alpha}}=\partial_{0}\varphi^{\bar{\alpha}}$,
see (\ref{eq:connTransLaw}). 
\end{rem}

\subsubsection{Covariant Derivative}

Given the connection (\ref{eq:Conn}) a covariant differential relative
to the connection $\Gamma$ can be introduced as a map $\nabla^{\Gamma}:\mathcal{J}^{1}\mathcal{(E)\rightarrow T}^{\ast}(\mathcal{B)\otimes V(E)}$
which reads in coordinates as
\begin{equation}
\nabla^{\Gamma}=\left(x_{0}^{\alpha}-\Gamma_{0}^{\alpha}\right)\mathrm{d}t^{0}\otimes\partial_{\alpha}.\label{eq:CovDiffDef}
\end{equation}
If the connection is trivial, i.e. $\Gamma_{0}^{\alpha}=0$ then the
covariant differential on the bundle $\mathcal{E}\rightarrow\mathcal{B}$
corresponds to the classical time derivative. Based on the covariant
differential one can define the covariant derivative of a section
$s:\mathcal{B\rightarrow E}$ which follows as 
\[
\nabla^{\Gamma}(s)=\left(\partial_{0}s^{\alpha}-\Gamma_{0}^{\alpha}\circ s\right)\mathrm{d}t^{0}\otimes\partial_{\alpha}.
\]
Given the vector field $\partial_{0}:\mathcal{B\rightarrow T(B)}$,
the contraction 
\[
\partial_{0}\rfloor\nabla^{\Gamma}(s)=\left(\partial_{0}s^{\alpha}-\Gamma_{0}^{\alpha}\circ s\right)\partial_{\alpha}
\]
is said to be the covariant derivative of $s:\mathcal{B\rightarrow E}$
along $\partial_{0}$, see \cite{Giachetta}.

\subsection{System Representation}

We define a time-variant (Port) Hamiltonian system on the bundle $\mathcal{Z}\rightarrow\mathcal{E}\rightarrow\mathcal{B}$
together with a connection (\ref{eq:Conn}) such that the structure
of the system is invariant with respect to bundle morphisms of $\mathcal{E}\rightarrow\mathcal{B}$
that explicitly depend on the time coordinate but we do not consider
time-reparametrization (i.e., $\bar{x}=\varphi(t,x)$ and $\bar{t}=t$).
The input bundle $\mathcal{U}\rightarrow\mathcal{E}$, $(t^{0},x^{\alpha},u^{i})\rightarrow(t^{0},x^{\alpha})$
is a vector bundle, and $\mathcal{Z}\rightarrow\mathcal{E}$ is used
in the affine case, analogous to the time-invariant case. 

Applying a time-variant transformation does not destroy the structure
of a (Port) Hamiltonian system if it is defined in an intrinsic manner.
But in contrast to \cite{FujimotoCSL}, where only the relative motion
with respect to a frame should be expressed in a Port Hamiltonian
framework, it is necessary in general to adopt the involved differential
operators to preserve the structure of the system. This will be done
by applying covariant derivatives. The choice of $\Gamma$ corresponds
to the selection of a frame of reference, i.e. starting in an inertial
frame (with trivial connection) a time-variant transformation induces
a non-trivial connection. An example how a non-trivial connection
can arise by changing the frame of reference for a mechanical system
is given in Section \ref{sub:Equations-of-motion}.
\begin{defn}
\label{def:Given-a-connection}Given a connection (\ref{eq:Conn})
then a Hamiltonian system on a bundle $\mathcal{Z}\rightarrow\mathcal{E}\rightarrow\mathcal{B}$
is given by
\begin{equation}
\partial_{0}\rfloor\nabla^{\Gamma}=\left(J-R\right)\rfloor\mathrm{d}H+G\rfloor u\label{eq:PCHTimeVarIntrinsic}
\end{equation}
with $J,R:\mathcal{V}^{\ast}\mathcal{(E)\rightarrow V(E)}$ where
$J$ is skew symmetric and $R$ is symmetric and positive semidefinite.
Additionally we have the bundle map $G:\mathcal{U}\rightarrow\mathcal{V(E)}$
and the Hamiltonian $H\in\mathcal{C}^{\infty}(\mathcal{E})$. In local
coordinates we obtain the expression 
\[
x_{0}^{\alpha}-\Gamma_{0}^{\alpha}=\left(J^{\alpha\beta}-R^{\alpha\beta}\right)\partial_{\beta}H+G_{i}^{\alpha}u^{i}.
\]
where 
\begin{equation}
v_{H,\mathcal{V}}^{\alpha}=\left(J^{\alpha\beta}-R^{\alpha\beta}\right)\partial_{\beta}H+G_{i}^{\alpha}u^{i}\label{eq:HamFieldRHS}
\end{equation}
is referred as the vertical part of the Hamiltonian vector field $v_{H,\mathcal{V}}=v_{H,\mathcal{V}}^{\alpha}\partial_{\alpha}$
with $H,G_{i}^{\alpha},J^{\alpha\beta},R^{\alpha\beta}\in\mathcal{C}^{\infty}\mathcal{(E)}$.
\end{defn}
To show that in the time-variant scenario $J$ and $R$ are maps from
the vertical cotangent bundle to the vertical tangent bundle $J,\, R:\mathcal{V}^{\ast}\mathcal{(E)\rightarrow V(E)}$
we have to observe that the total differential of the Hamiltonian
$\mathrm{d}H=\partial_{0}H\mathrm{d}t^{0}+\partial_{\beta}H\mathrm{d}x^{\beta}$
can be decomposed according to (\ref{eq:splitvec1a}) as 
\begin{equation}
\mathrm{d}H=\partial_{\beta}H\omega_{\mathcal{V}}^{\beta}+w_{0}^{\mathcal{H}}(H)\mathrm{d}t^{0}\label{eq:dHhorandverdeco}
\end{equation}
with the \textit{1-form} $\omega_{\mathcal{V}}^{\alpha}$%
\footnote{It should be noted that $\omega_{\mathcal{V}}^{\alpha}=\mathrm{d}x^{\alpha}-\Gamma_{0}^{\alpha}\mathrm{d}t^{0}$
serves as an adapted basis for $\mathcal{V}^{*}(\mathcal{E})$ induced
by the connection $\Gamma$. This justifies $J,\, R:\mathcal{V}^{\ast}\mathcal{(E)\rightarrow V(E)}$
as we have claimed. %
} and the vector field $w_{0}^{\mathcal{H}}$ from (\ref{eq:whwv})
from which we deduce the result, together with (\ref{eq:CovDiffDef}). 
\begin{rem}
Since the exterior derivative $\mathrm{d}$ can also be decomposed
into a vertical and a horizontal one $\mathrm{d}=\mathrm{d}_{\mathcal{V}}+\mathrm{d}_{\mathcal{H}}$,
see \cite{Giachetta,Saunders} one can equivalently use $\mathrm{d}_{\mathcal{V}}H$
$ $instead of $\mathrm{d}H$ in the relation (\ref{eq:PCHTimeVarIntrinsic}).
\end{rem}
Let us consider the total time change $d_{0}(H)$ of the Hamiltonian
$H$ along solutions. For a time-invariant system we obtain a relation
of the form
\[
d_{0}(H)=v_{H}(H)=-(\partial_{\alpha}H)R^{\alpha\beta}(\partial_{\beta}H)+u^{i}y_{i}
\]
where $v_{H}(H)$ denotes the Lie-derivative of the Hamiltonian $H$
along the Hamiltonian vector field $v_{H}=\left(J^{\alpha\beta}-R^{\alpha\beta}\right)\partial_{\beta}H+G_{i}^{\alpha}u^{i}$,
see (\ref{eq:HamTimeInvStandard}) and the total time derivative (time-invariant
case) reads as $d_{0}=\dot{x}^{\alpha}\partial_{\alpha}$.\\
The intrinsic version in the time-variant case is given in the following
corollary. 
\begin{cor}
The total time change of the Hamiltonian leads to a decomposition
of the form 
\begin{equation}
d_{0}(H)=w_{0}^{\mathcal{H}}(H)+v_{H,\mathcal{V}}\left(H\right)\label{eq:DecompHdot}
\end{equation}
or in local coordinates
\[
\underset{w_{0}^{\mathcal{H}}(H)}{\underbrace{\partial_{0}(H)+\Gamma_{0}^{\alpha}\partial_{\alpha}(H)}}\underset{v_{H,\mathcal{V}}\left(H\right)}{\underbrace{-(\partial_{\alpha}H)R^{\alpha\beta}(\partial_{\beta}H)+\partial_{\alpha}(H)G_{i}^{\alpha}u^{i}}}
\]
where the vector field $w_{0}^{\mathcal{H}}$ is induced by the connection
$\Gamma$, see (\ref{eq:whwv}), the vector field $v_{H,\mathcal{V}}$
from (\ref{eq:HamFieldRHS}) and $d_{0}=\partial_{0}+x_{0}^{\alpha}\partial_{\alpha}$.
It is worth noting that we have a decomposition of the derivative
operator into a horizontal $w_{0}^{\mathcal{H}}$ and a vertical $v_{H,\mathcal{V}}$
component where the horizontal part degenerates to $\partial_{0}$
only if the connection is trivial, i.e. $\Gamma_{0}^{\alpha}=0$.
\end{cor}
This follows by a direct calculation based on the intrinsic system
representation. The discussion of collocation will be performed using
a special class of systems in the next section.

\subsection{A special class of Port Hamiltonian Systems}

Mechanical systems in a Port Hamiltonian representation are distinguished
since the bundle $\mathcal{E}\rightarrow\mathcal{B}$ has an even
richer geometric structure. The main difference is the separation
of the $x^{\alpha}$ coordinates in positions $q^{\alpha}$ and momenta
$\dot{q}_{\alpha}$. To obtain the correct geometric picture we introduce
the bundle $\mathcal{Q}\rightarrow\mathcal{B}$ with coordinates $(t^{0},q^{\alpha})$
for $\mathcal{Q}$. This bundle structure induces again some tangent
structures where we make use of the dual vertical cotangent bundle
$\mathcal{V}^{\ast}\mathcal{(Q)\rightarrow Q}$ with coordinates $(t^{0},q^{\alpha},\dot{q}_{\alpha})$.
Furthermore we will utilize $\mathcal{T}(\mathcal{V}^{\ast}(\mathcal{Q}))$
which possesses the adapted bases $(\partial_{0},\partial_{\alpha},\dot{\partial}^{\alpha})$
with $\dot{\partial}^{\alpha}=\frac{\partial}{\partial\dot{q}_{\alpha}}$
and to be conform with most of the literature we make the identification
$\dot{q}_{\alpha}=p_{\alpha}$.

\subsubsection{The Composite Bundle Structure}

From the bundles introduced so far $\mathcal{V}^{\ast}\mathcal{(Q)\rightarrow Q}$
and $\mathcal{Q\rightarrow B}$ we can construct the composite bundle
structure $\mathcal{V}^{\ast}\mathcal{(Q)\rightarrow Q}\rightarrow\mathcal{B}$. 

The additional fibration $\mathcal{V}^{\ast}\mathcal{(Q)\rightarrow B},\,(t^{0},q^{\alpha},p_{\alpha})\rightarrow(t^{0})$
plays the role of the state bundle $\mathcal{E}\rightarrow\mathcal{B}$.
We choose a connection $\gamma$ that corresponds to the selection
of a frame of reference that splits $\mathcal{T(Q)}$ 
\begin{equation}
\gamma=\mathrm{d}t^{0}\otimes\left(\partial_{0}+\gamma_{0}^{\alpha}\partial_{\alpha}\right)\,,\,\gamma_{0}^{\alpha}\in C^{\infty}(\mathcal{Q}).\label{eq:Mech Space Conn}
\end{equation}
Based on the connection $\gamma$ the so-called covertical connection,
see \cite{Giachetta} 
\begin{equation}
\Gamma_{H}=\mathrm{d}t^{0}\otimes\left(\partial_{0}+\gamma_{0}^{\alpha}\partial_{\alpha}-(\partial_{\rho}\gamma_{0}^{\beta})p_{\beta}\dot{\partial}^{\rho}\right)\ \label{Ham connection}
\end{equation}
can be constructed that splits the tangent bundle $\mathcal{T}(\mathcal{V}^{\ast}(\mathcal{Q}))$
with respect to the fibration $\mathcal{V}^{\ast}(\mathcal{Q})\rightarrow\mathcal{B}$
and this is the connection that corresponds to $\Gamma$, see (\ref{eq:Conn})
that splits $ $$\mathcal{T}(\mathcal{E})$. \\
In classical mechanics it is well known that the Hamiltonian vector
field can be defined using Symplectic and Poisson structures, see
\cite{Abraham,Giachetta}. These concepts can also be generalized
to the time-variant scenario. We base the construction of the Hamiltonian
vector field using the Hamilton form $\omega_{H}$ which follows from
the canonical Liouville \textit{1-form }(by a pull-back)\textit{ }
\[
\omega_{H}=p_{\alpha}\mathrm{d}q^{\alpha}-(H+p_{\alpha}\gamma_{0}^{\alpha})\mathrm{dt^{0}}
\]
with $H\in C^{\infty}(\mathcal{V}^{\ast}\mathcal{(Q)})$. \\
From the relations $v_{H}\rfloor\mathrm{d}\omega_{H}=0$ and $v_{H}\rfloor\mathrm{d}t^{0}=1$
the autonomous Hamiltonian vector field $v_{H}:\mathcal{V}^{\ast}\mathcal{(Q)}\rightarrow\mathcal{T(V}^{\ast}\mathcal{(Q)}$
follows with $H_{\gamma}=H+p_{\beta}\gamma_{0}^{\beta}$ as 
\begin{equation}
v_{H}=\partial_{0}+\dot{\partial}^{\alpha}H_{\gamma}\partial_{\alpha}-\partial_{\alpha}H_{\gamma}\dot{\partial}^{\alpha}.\label{Ham vector field}
\end{equation}

\begin{rem}
If the connection $\gamma$ is trivial, i.e. $\gamma_{0}^{\alpha}=0$
holds then the Hamiltonian vector field reads as $v_{H}=\partial_{0}+\dot{\partial}^{\alpha}H\partial_{\alpha}-\partial_{\alpha}H\dot{\partial}^{\alpha}$
since then $H_{\gamma}=H$ which is the standard result in mechanics,
see \cite{Abraham}.
\end{rem}

\subsubsection{Power Balance Equation for Controlled Mechanical Systems}

The connection (\ref{Ham connection}) enables us to split the Hamiltonian
vector field (\ref{Ham vector field}) into a vertical and horizontal
part, respectively, see also \cite{SchoeberlPAMM2006}. It follows
that (analogously to (\ref{eq:splitvec1a})) we have a decomposition
of the form
\begin{equation}
\begin{array}{ccl}
v_{H,\mathcal{V}} & = & \dot{\partial}^{\alpha}H\partial_{\alpha}-\partial_{\alpha}H\dot{\partial}^{\alpha}\\
v_{H,\mathcal{H}} & = & \partial_{0}+\gamma_{0}^{\alpha}\partial_{\alpha}-(\partial_{\rho}\gamma_{0}^{\beta})p_{\beta}\dot{\partial}^{\rho}
\end{array}\label{eq:ham split}
\end{equation}
according to the Hamiltonian connection (\ref{Ham connection}). $ $Let
us inspect the relation (\ref{eq:DecompHdot}) in the case of a mechanical
system where it is obvious that now the vector field $ $$w_{0}^{\mathcal{H}}$
from (\ref{eq:whwv}) is replaced by $v_{H,\mathcal{H}}$ due to our
richer bundle structure. We consider the autonomous case first and
since no dissipation is present we obtain
\begin{eqnarray}
v_{H}(H) & = & v_{H,\mathcal{H}}(H)\label{eq:Ham time change hor}
\end{eqnarray}
where $v_{H}(H)$ denotes the Lie-derivative of the \textit{Hamiltonian
}$H$ with respect to the vector field $v_{H}$. To include control
inputs we consider an extended Hamiltonian corresponding to a controlled
Hamiltonian, see \cite{NijmeijervanderSchaft}, of the form
\begin{equation}
H=H_{0}-H_{c,\rho}u^{\rho}\label{ham extended input}
\end{equation}
with $H_{0},\, H_{c,\rho}\in C^{\infty}(\mathcal{V}^{\ast}\mathcal{(Q)})$
and the input functions $u^{\rho}\in\mathcal{C}^{\infty}(\mathcal{B})$. 

Based on (\ref{ham extended input}) and the relation (\ref{eq:Ham time change hor})
we can state the following Theorem concerning a covariant formulation
of the power flows as well as the concept of collocation. 
\begin{thm}
\label{thm:The-change-of}The change of the free Hamiltonian $H_{0}$
along solutions of the Hamiltonian system decomposes as
\begin{eqnarray*}
v_{H}(H_{0}) & = & v_{H,\mathcal{H}}(H_{0})+v_{H,\mathcal{V}}(H_{c,\rho})u^{\rho}
\end{eqnarray*}
where we used the controlled Hamiltonian $H=H_{0}-H_{c,\rho}u^{\rho}$
and the decomposition (\ref{eq:ham split}).
\end{thm}
Obviously the choice of the output $y_{\rho}=v_{H,\mathcal{V}}(H_{c,\rho})$
allows a physical interpretation of the power flows of the system,
since $v_{H,\mathcal{H}}(H_{0})$ corresponds to the power caused
by the free Hamiltonian $H_{0}$ and the product $y_{\rho}u^{\rho}$
describes the power flow into the system caused by the input. Theorem
\ref{thm:The-change-of} is an intrinsic version of the balancing/interaction
of power flows and it reduces to the well known formula when the connection
is trivial, i.e. $\gamma_{0}^{\alpha}=0$. It is worth mentioning
that the splitting of the Hamiltonian field as in (\ref{eq:ham split})
is essential to obtain a coordinate free representation.

\section{Applications}

This section is devoted to a short discussion of two possible applications
of the presented theory. Firstly, we focus on the description of the
error system arising when the stabilization of the trajectory tracking
error is the objective in a Port Hamiltonian framework and discuss
the role of the connection in this context. As a second application
we present the equations of motion of a mass particle observed from
a rotating frame of reference, in a time-variant Hamiltonian formulation
since this example demonstrates how to calculate the non-trivial connections
coefficients and additionally one can show that this connection is
used to formulate conservations laws and/or power balance relations.

\subsection{Trajectory tracking\label{sub:Trajectory-tracking}}

This section deals with time-variant Port Hamiltonian systems as they
arise quite naturally when a feed-forward based approach is applied
to a time-independent Port Hamiltonian system and the control objective
is the stabilization of the error system using techniques of passivity.
It is evident that the error system is constructed using a time-variant
transformation which enables us to discuss the developed machinery
on this concrete example. It should be pointed out that compared to
the approach in \cite{FujimotoAut} we use a different definition
of a time-variant Port Hamiltonian system (see Definition \ref{def:Given-a-connection})
which is coordinate independent and more general due to the connection
term. If the connection is interpreted as an additive term in a modified
Hamiltonian then our results coincide with the results of \cite{FujimotoAut}.

Let us consider a Port Hamiltonian system of the form 
\begin{equation}
x_{0}^{\alpha}=\left(J^{\alpha\beta}-R^{\alpha\beta}\right)\partial_{\beta}H+G_{i}^{\alpha}u^{i}.\label{eq:controlHamSysOrig}
\end{equation}
The control system (\ref{eq:controlHamSysOrig}) is modeled on a bundle
$\mathcal{E}\rightarrow\mathcal{B}$ where the coordinates are obviously
adapted to the connection, such that the connection coefficients then
read as $\Gamma_{0}^{\alpha}=0.$ If a desired trajectory $c_{d}(t)$
and a corresponding input $\eta_{d}(t)$, which produces this trajectory,
is given 
\begin{equation}
\partial_{0}c_{d}^{\alpha}=\left(\left(J^{\alpha\beta}-R^{\alpha\beta}\right)\partial_{\beta}H+G_{i}^{\alpha}\eta_{d}^{i}\right)\circ c_{d}\label{eq:controlHamSysOrigSol}
\end{equation}
then a transformation with respect to a reference trajectory can be
stated as 
\begin{eqnarray}
\bar{x}^{\bar{\alpha}} & = & \varphi^{\bar{\alpha}}(t^{0},x^{\beta})=\delta_{\alpha}^{\bar{\alpha}}(x^{\alpha}-c_{d}^{\alpha})\label{eq:controlxtrans}\\
\bar{u}^{\bar{j}} & = & M_{i}^{\bar{j}}(u^{i}-\eta_{d}^{i})\,,\,\,\,\, M_{i}^{\bar{j}}\in C^{\infty}(\mathcal{X}).\label{eq:controlutrans}
\end{eqnarray}

\begin{rem}
The relations (\ref{eq:controlxtrans}) and (\ref{eq:controlutrans})
are of course only a trivial choice leading to an error system. But
for our purposes, the geometric interpretation of time-variant Port
Hamiltonian systems, the exact structure of (\ref{eq:controlxtrans})
and (\ref{eq:controlutrans}) is of minor importance, since the intrinsic
system formulation is independent of the special choice of the transition
functions. For a more general discussion of error systems in the Port
Hamiltonian context, see \cite{FujimotoAut}. 
\end{rem}
It is readily observed that the connection coefficients read as 
\begin{equation}
\bar{\Gamma}_{0}^{\bar{\alpha}}=\partial_{0}\varphi^{\bar{\alpha}}=-(\partial_{0}c_{d}^{\alpha})\delta_{\alpha}^{\bar{\alpha}}.\label{eq:controlConnectionCoeff}
\end{equation}
with $\varphi^{\bar{\alpha}}$ from (\ref{eq:controlxtrans}). \\
The combination of the results of Lemma \ref{lem:affineInput} together
with the fact that in the time-variant case a connection appears that
leads to an additional affine term, we obtain the following Corollary. 
\begin{cor}
\label{cor:Given--and}Given $c_{d}(t)$ and $\eta_{d}(t)$ fulfilling
(\ref{eq:controlHamSysOrigSol}) then the transformations (\ref{eq:controlxtrans})
and (\ref{eq:controlutrans}) applied to the system (\ref{eq:controlHamSysOrig})
lead to a representation as
\begin{equation}
\bar{x}_{0}^{\bar{\alpha}}=(\bar{J}^{\bar{\alpha}\bar{\beta}}-\bar{R}^{\bar{\alpha}\bar{\beta}})\partial_{\bar{\beta}}(\bar{H}+\breve{H})+\bar{G}_{\bar{i}}^{\bar{\alpha}}\bar{u}^{\bar{i}},\label{eq:controlHamTransSys}
\end{equation}
if and only if the partial differential equations 
\[
(\bar{J}^{\bar{\alpha}\bar{\beta}}-\bar{R}^{\bar{\alpha}\bar{\beta}})\partial_{\bar{\beta}}\breve{H}=\underset{*}{\underbrace{\bar{\Gamma}_{0}^{\bar{\alpha}}}}+\underset{**}{\underbrace{\delta_{\alpha}^{\bar{\alpha}}G_{i}^{\alpha}\eta_{d}^{i}}}
\]
allow a solution for $\breve{H}$. 
\end{cor}
From Corollary \ref{cor:Given--and} we deduce that if the connection
{*} and the feed-forward part {*}{*} are expressed as an additive
Hamiltonian, one obtains a Port-Hamiltonian representation which is
beneficial when the control objective is to stabilize error systems
that arise typically when trajectory tracking is the demand

\subsection{Rotating Frame of Reference\label{sub:Equations-of-motion}}

To show that a time-variant transformation preserves the Hamiltonian
structure we will consider the equations of motion of a mass particle
using two different frames, i.e. an inertial one and a rotating one
with respect to the inertial frame. Let us consider an inertial system
with Euclidean coordinates $(q^{\alpha},t^{0})$ together with the
bundle structure $\mathcal{Q}\rightarrow\mathcal{B}$. The canonical
equations of motion for a mass particle with mass $m\in\mathbb{R}^{+}$
read as 
\begin{eqnarray}
\partial_{0}s^{\alpha}=\dot{\partial}^{\alpha}H & \,,\, & \partial_{0}(p_{\alpha}\circ s)=-\partial_{\alpha}H\label{eq:HamMassInert}
\end{eqnarray}
with $s:\mathcal{B}\rightarrow\mathcal{Q}$. For this example the
Hamiltonian is given as 
\[
H=\frac{1}{2m}p_{\alpha}\delta^{\alpha\beta}p_{\beta}
\]
where $\delta$ denotes the Kronecker delta. \\
A rotating coordinate chart with respect to the inertial one can be
constructed, using $\bar{q}^{\bar{\alpha}}=R_{\beta}^{\bar{\alpha}}(t^{0})q^{\beta}=\varphi^{\bar{\alpha}}(q^{\beta},t^{0})$
with $R_{\beta}^{\bar{\alpha}}\delta^{\alpha\beta}R_{\alpha}^{\bar{\beta}}=\delta^{\bar{\alpha}\bar{\beta}}$.
The non-trivial connection in the floating reference system can be
computed as in (\ref{eq:connTransLaw}) and reads as $\bar{\gamma}_{0}^{\bar{\alpha}}=\partial_{0}(R_{\beta}^{\bar{\alpha}})R_{\bar{\rho}}^{\beta}\bar{q}^{\bar{\rho}}=\Omega_{\bar{\rho}}^{\bar{\alpha}}\bar{q}^{\bar{\rho}}.$ 

The equations of motion in the rotating coordinate system follow as,
see (\ref{Ham vector field}) or \cite{SchoeberlJMathP} for a detailed
exposition concerning covariant derivatives
\begin{equation}
\begin{array}{rcc}
\partial_{0}\bar{s}^{\bar{\alpha}}-\bar{\gamma}_{0}^{\bar{\alpha}} & = & \dot{\partial}^{\bar{\alpha}}\bar{H}\\
\partial_{0}(\bar{p}_{\bar{\alpha}}\circ\bar{s})+\bar{p}_{\bar{\beta}}\partial_{\bar{\alpha}}\bar{\gamma}_{0}^{\bar{\beta}} & = & -\partial_{\bar{\alpha}}\bar{H}.
\end{array}\label{eq:HamMassRot-1}
\end{equation}
Comparing (\ref{eq:HamMassInert}) with (\ref{eq:HamMassRot-1}) it
is evident that they both are Hamiltonian representations of the same
physical problem, but in contrast to the inertial frame where the
left hand side consists of partial time derivatives only in the rotating
frame a differential operator induced by the connection (\ref{Ham connection})
has to be applied. In the concrete example where the connection reads
as $\bar{\gamma}_{0}^{\bar{\alpha}}=\Omega_{\bar{\rho}}^{\bar{\alpha}}\bar{q}^{\bar{\rho}}$
together with $\bar{s}:\mathcal{B}\rightarrow\bar{\mathcal{Q}}$ the
equations follow as
\begin{eqnarray}
\partial_{0}\bar{s}^{\bar{\alpha}}-\Omega_{\bar{\rho}}^{\bar{\alpha}}\bar{s}^{\bar{\rho}} & = & \frac{1}{m}\bar{p}_{\bar{\beta}}\delta^{\bar{\alpha}\bar{\beta}}\label{eq:HamRobEq1}\\
\partial_{0}(\bar{p}_{\bar{\alpha}}\circ\bar{s})+\bar{p}_{\bar{\beta}}\Omega_{\bar{\alpha}}^{\bar{\beta}} & = & 0,\label{eq:HamRobEq2}
\end{eqnarray}
when no other forces are applied.
\begin{rem}
To apply Theorem \ref{thm:The-change-of} it should be noted that
$\bar{v}_{\bar{H},\mathcal{H}}$ has to be constructed with the connection
$\bar{\gamma}_{0}^{\bar{\alpha}}=\Omega_{\bar{\rho}}^{\bar{\alpha}}\bar{q}^{\bar{\rho}}.$
Then the Hamiltonian $\bar{H}$ is a conserved quantity as the conservation
law $\bar{v}_{\bar{H},\mathcal{H}}(\bar{H})=0$ (see Theorem \ref{thm:The-change-of}
with $\bar{H}_{0}=\bar{H}$) $ $is met. If we apply forces to control
the system then the full relation of Theorem \ref{thm:The-change-of}
has to be applied including the additional expression $\bar{v}_{\bar{H},\mathcal{V}}$,
see (\ref{eq:ham split}). 
\end{rem}
Combining the equations (\ref{eq:HamRobEq1}) and (\ref{eq:HamRobEq2})
we obtain 
\[
\partial_{0}(\partial_{0}\bar{s}^{\bar{\gamma}}-\Omega_{\bar{\rho}}^{\bar{\gamma}}\bar{s}^{\bar{\rho}})\delta_{\bar{\gamma}\bar{\alpha}}+(\partial_{0}\bar{s}^{\bar{\gamma}}-\Omega_{\bar{\rho}}^{\bar{\gamma}}\bar{s}^{\bar{\rho}})\delta_{\bar{\gamma}\bar{\beta}}\Omega_{\bar{\alpha}}^{\bar{\beta}}=0
\]
which can be rewritten as 
\[
\partial_{00}\bar{s}^{\bar{\beta}}-\partial_{0}\Omega_{\bar{\rho}}^{\bar{\beta}}\bar{s}^{\bar{\rho}}-2\Omega_{\bar{\rho}}^{\bar{\beta}}\partial_{0}\bar{s}^{\bar{\rho}}+\Omega_{\bar{\gamma}}^{\bar{\beta}}\Omega_{\bar{\rho}}^{\bar{\gamma}}\bar{s}^{\bar{\rho}}=0
\]
where we have exploited the skew symmetry of the so-called angular
velocity tensor $\Omega$. This is a classical result in mechanics,
i.e. this is an expression for the acceleration including the Coriolis
and the Centrifugal acceleration. The relations (\ref{eq:HamRobEq1})
and (\ref{eq:HamRobEq2}) are the covariant version in a Hamiltonian
point of view. 

Summarizing we can state, that time-variant transformations preserve
the Hamiltonian structure, if a (Port) Hamiltonian system is introduced
covariantly as in Definition \ref{def:Given-a-connection}. The system
properties can be expressed in a time-variant frame which may lead
to non trivial connections, which then requires the use of covariant
derivatives. Also the conservation laws and the power balance laws
have to be formulated in an intrinsic fashion, see Theorem \ref{thm:The-change-of}.

Finally it should be stressed again, that the interpretation of the
connection in the subsections (\ref{sub:Trajectory-tracking}) and
(\ref{sub:Equations-of-motion}) is completely different. In subsection
(\ref{sub:Trajectory-tracking}) the connection is absorbed in a modified
Hamiltonian, whereas in subsection (\ref{sub:Equations-of-motion})
the connection is explicitly part of the covariant differential.

\section{Conclusion}

In this paper we investigated the geometry of time-variant (Port)
Hamiltonian systems and used an intrinsic description based on connections
and covariant derivatives. We were interested in two concrete applications
where time-variant systems arise quite naturally, namely error systems
with regard to trajectory tracking in control theory as well as time-variant
mechanics. Concerning control theory we described that the formulation
of the error system in a Hamiltonian fashion requires to absorb the
nontrivial connection and the feed-forward part in a modified Hamiltonian.
Here the key problem was the definition of a time-variant Hamiltonian
system itself, i.e. if the tensors are time-dependent and/or the coordinate
chart is moving. The interpretation of these two scenarios is significant
for a correct understanding of time-variant systems. Furthermore time-variant
mechanics has been analyzed in a covariant way, where we showed how
collocation and the balancing/interaction of power flows can be formulated
in an intrinsic way. 

\section*{Acknowledgment}

Markus Sch\"{o}berl is an APART fellowship holder of the Austrian Academy of Sciences and he was partially supported
by the Austrian center of competence in Mechatronics (ACCM). The authors want to thank the anonymous reviewers for their helpful 
comments and suggestions to improve the quality and the readability of the present paper.
\bibliographystyle{plain}   
\bibliography{./maxibib1}

\begin{thebibliography}{10}

\bibitem{Abraham}
R.A. Abraham and J.E. Marsden.
\newblock {\em Foundations of Mechanics}, volume 2nd ed.
\newblock Addison-Wesley, Reading, MA, 1978.

\bibitem{Cheng}
D.~Cheng, A.~Astolfi, and R.~Ortega.
\newblock On feedback equivalence to port controlled hamiltonian systems.
\newblock {\em Sys. Control Lett.}, 54:911--917, 2005.

\bibitem{FujimotoCSL}
K.~Fujimoto and T.~Sugie.
\newblock Canonical transformation of generalized hamiltonian systems.
\newblock {\em Sys. Control Lett.}, 42:217--227, 2001.

\bibitem{FujimotoAut}
K.~Fujimoto and T.~Sugie.
\newblock Trajectory tracking of port-controlled hamiltonian systems via
  generalized canonical transformations.
\newblock {\em Automatica}, 39:2059--2069, 2003.

\bibitem{Giachetta}
G.~Giachetta, G.~Sardanashvily, and L.~Mangiarotti.
\newblock {\em New Lagrangian and Hamiltonian Methods in Field Theory}.
\newblock World Scientific, 1997.

\bibitem{Gotay}
M.J. Gotay.
\newblock A multisymplectic framework for classical field theory and the
  calculus of variations ii: Space + time decomposition.
\newblock {\em Differential Geometry and its Applications}, (1):375--390, 1991.

\bibitem{Kanat}
I.V. Kanatchikov.
\newblock Canonical structure of classical field theory in the polymomentum
  phase space.
\newblock {\em Rep. on Math. Physics}, 41(1):49--90, 1998.

\bibitem{MaschkeOrtegaSchaft2000TAC}
B.~Maschke, R.~Ortega, and A.J. {van der Schaft}.
\newblock Energy-based lyapunov functions for forced hamiltonian systems with
  dissipation.
\newblock {\em IEEE Trans. on Aut. Contr.}, 45:1498--1502, 2000.

\bibitem{NijmeijervanderSchaft}
H.~Nijmeijer and A.J. {van der Schaft}.
\newblock {\em Nonlinear Dynamical Control Systems}.
\newblock Springer-Verlag, New York, 1990.

\bibitem{OrtegaEnergy}
R.~Ortega, A.J. {van der Schaft}, I.~Mareels, and B.~Maschke.
\newblock Putting energy back in control.
\newblock {\em IEEE Contr. Syst. Mag}, 21:18--33, 2001.

\bibitem{Ortega}
R.~Ortega, A.J. {van der Schaft}, B.~Maschke, and G.~Escobar.
\newblock Interconnection and damping assignment passivity-based control of
  port hamiltonian systems.
\newblock {\em Automatica}, 36:585--596, 2002.

\bibitem{Saunders}
D.~J. Saunders.
\newblock {\em The Geometry of Jet Bundles}.
\newblock Mathematical University Press, Cambridge University Press, Cambridge,
  1989.

\bibitem{SchoeberlPAMM2006}
M.~Sch{\"o}berl and K.~Schlacher.
\newblock Geometric analysis of hamiltonian mechanics using connections.
\newblock In {\em PAMM, Proceedings of GAMM}, pages 843--844, 2006.

\bibitem{SchoeberlJMathP}
M.~Sch{\"o}berl and K.~Schlacher.
\newblock Covariant formulation of the governing equations of continuum
  mechanics in an eulerian description.
\newblock {\em J. Math. Phys}, 48(5), 2007.

\bibitem{SchoeberlNolcos}
M.~Sch{\"o}berl, R.~Stadlmayr, and K.~Schlacher.
\newblock Geometric analysis of time variant hamiltonian control systems.
\newblock In {\em Proceedings IFAC Symposium on Nonlinear Control Systems,
  NOLCOS 2007}, pages 1026--1031, 2007.

\bibitem{vanderSchaft}
A.J. {van der Schaft}.
\newblock {\em L2-Gain and Passivity Techniques in Nonlinear Control}.
\newblock Springer-Verlag, New York, 2000.

\end{thebibliography}

\end{document}